\documentclass[a4paper]{amsart}

\usepackage{amsmath}
\usepackage{amssymb}
\usepackage{bm} 
\usepackage{mathrsfs} 
\usepackage{tensor} 
\usepackage{braket} 
\usepackage{physics} 

\DeclareMathOperator{\Vol}{Vol}

\DeclareMathOperator{\Ric}{Ric}

\DeclareMathOperator{\End}{End}

\DeclareMathOperator{\Exp}{Exp}
\DeclareMathOperator{\Area}{Area}

\newcommand{\clop}[2]{[#1, #2)}

\usepackage{amsthm} 

\usepackage{cleveref} 

\usepackage{autonum} 

\crefname{theorem}{Theorem}{Theorems}
\crefname{lemma}{Lemma}{Lemmas}
\crefname{proposition}{Proposition}{Propositions}
\crefname{corollary}{Corollary}{Corollaries}
\crefname{claim}{Claim}{Claims}
\crefname{assumption}{Assumption}{Assumptions}
\crefname{definition}{Definition}{Definitions}
\crefname{remark}{Remark}{Remarks}
\crefname{example}{Example}{Examples}
\crefname{equation}{}{}
\crefname{section}{Section}{Sections}
\crefname{subsection}{Section}{Sections}
\crefname{appendix}{Appendix}{Appendices}

\crefname{condition}{}{}

\creflabelformat{equation}{(#2#1#3)}
\creflabelformat{condition}{(#2#1#3)}

\theoremstyle{plain}
\newtheorem{theorem}{Theorem}[section]

\newtheorem{proposition}[theorem]{Proposition}

\theoremstyle{definition}

\theoremstyle{remark}

\numberwithin{equation}{section}

\usepackage{comment}
\usepackage{tikz-cd}

\newcommand{\RN}[1]{%
	\textup{\uppercase\expandafter{\romannumeral#1}}%
}

\title{On the second variation of the Graham-Witten energy}
\author{Yuya Takeuchi}
\address{Graduate School of Mathematical Sciences \\ The University of Tokyo \\ 3-8-1 Komaba, Meguro, Tokyo 153-8914 Japan}
\curraddr{Department of Mathematics \\ Graduate School of Science \\ Osaka University \\ 1-1 Machikaneyama-cho, Toyonaka, Osaka 560-0043, Japan}
\email{yu-takeuchi@cr.math.sci.osaka-u.ac.jp}



\begin{document}

\begin{abstract}
	The area renormalization procedure gives
	an invariant of even-dimensional closed submanifolds in a conformal manifold,
	which we call the Graham-Witten energy,
	and it is a generalization of the classical Willmore energy.
	In this paper,
	we obtain an explicit formula for the second variation of this energy
	at minimal submanifolds in an Einstein manifold.
	As an application,
	we prove that the even-dimensional totally geodesic spheres in the unit sphere
	are critical points of the Graham-Witten energy with non-negative second variation.
\end{abstract}

\maketitle

\section{Introduction} \label{section:introduction}

The volume renormalization for asymptotically hyperbolic Einstein manifolds
has played a crucial role in conformal geometry
after introduced by Henningson and Skenderis~\cite{Henningson-Skenderis98}
in studies of the AdS/CFT correspondence.
Let $\overline{X}$ be an $(n + 1)$-dimensional manifold with smooth boundary $M$,
and $r$ be a defining function of the boundary.
A Riemannian metric $g^{+}$ on the interior $X$ of $\overline{X}$
is called an asymptotically hyperbolic Einstein metric
if $r^{2} g^{+}$ extends to a smooth Riemannian metric on $\overline{X}$
and $\Ric(g^{+}) + n g^{+} = O(r^{n - 2})$ holds.
The conformal class $C = [(r^{2} g^{+})|_{M}]$ on $M$
is independent of the choice of $r$,
called the conformal infinity of $g^{+}$.
If we choose a representative $g$ of $C$,
there exist a unique defining function $r$ and an identification of a neighborhood of $M$
with $\clop{0}{\epsilon_{0}} \times M$ such that
$g^{+}$ is of the form
\begin{equation}
	g^{+}
	= \frac{d r^{2} + g_{r}}{r^{2}},
\end{equation}
where $g_{r}$ is a one-parameter family of metrics on $M$ with $g_{0} = g$.
For such an $r$,
the volume $\Vol(\Set{ r > \epsilon} )$ of the domain $\Set{ r > \epsilon }$
has the following expansion, as $\epsilon \to + 0$:
\begin{equation} \label{eq:volume-renormalization}
	\Vol(\Set{ r > \epsilon })
	= \sum_{l = 0}^{\lceil n/2 \rceil -1} a_{l} \, \epsilon^{- n + 2 l}
	+
	\begin{cases}
		V^{o} + o(1), & n \ \text{odd}, \\
		L \log (1 / \epsilon) + V^{e} + o(1), & n \ \text{even}.
	\end{cases}
\end{equation}
The constant term $V^{o}$ or $V^{e}$ in this asymptotics is called the renormalized volume.
Graham~\cite[Theorem 3.1]{Graham00}
has proved that $V^{o}$ and $L$ are independent of the choice of $g$.
Moreover,
Graham and Zworski~\cite[Theorem 3]{Graham-Zworski03} have shown
that
\begin{equation} \label{eq:logarithmic-coefficient-of-volume-renormalization}
	 L = 2 c_{n / 2} \overline{Q},
\end{equation}
where $c_{k} = (- 1)^{k} [2^{2 k} (k -1)! k!]^{-1}$
and $\overline{Q}$ is the the total $Q$-curvature,
a conformal invariant of the boundary;
see also~\cite[Theorem 3.2]{Fefferman-Graham02}.

Graham and Witten~\cite{Graham-Witten99}
have introduced the \emph{area renormalization},
a renormalization procedure concerning submanifolds,
or more generally immersions.
It has been of great importance in physics recently
since it is closely related to the entanglement entropy in conformal field theory~\cite{Ryu-Takayanagi06}.
To simplify the exposition,
we consider only a closed $d$-dimensional submanifold $\Sigma$ in $M$ here;
see \cref{section:Graham-Witten-invariant} for the case of immersions.
There exists a $(d + 1)$-dimensional ``asymptotically minimal'' submanifold $Y$ in $X$
satisfying $Y \cap M = \Sigma$.
The area $\Area(Y \cap \Set{ r > \epsilon })$ of $Y \cap \Set{ r > \epsilon }$ with respect to $g^{+}$
has a similar asymptotic expansion to \cref{eq:volume-renormalization},
as $\epsilon \to + 0$:
\begin{equation}
	\Area(Y \cap \Set{ r > \epsilon })
	= \sum_{l = 0}^{\lceil d/2 \rceil -1} b_{l} \, \epsilon^{- d + 2 l}
	+
	\begin{cases}
		A^{o} + o(1), & d \ \text{odd}, \\
		K \log (1 / \epsilon) + A^{e} + o(1), & d \ \text{even}.
	\end{cases}
\end{equation}
The constant term $A^{o}$ or $A^{e}$ is called the renormalized area.
Graham and Witten~\cite[Proposition 2.1]{Graham-Witten99}
have proved that $A^{o}$ and $K$ are independent of the choice of the representative $g$.
In this paper,
the \emph{Graham-Witten energy} $\mathcal{E}$ of an even-dimensional closed submanifold $\Sigma$
is defined as
\begin{equation}
	K
	= 2 c_{d / 2} \mathcal{E}.
\end{equation}
This is a conformal invariant of the submanifold $\Sigma$,
and coincides with the classical Willmore energy
in the case of two-dimensional submanifolds in $\mathbb{R}^{3}$ or $S^{3}$.
Note that
our definition of the Graham-Witten energy is different from
that in \cite{Graham-Reichert} and \cite{Zhang} by a constant factor.

The first variation of the Graham-Witten energy has been recently obtained by
Graham and Reichert~\cite[Theorem 4.3]{Graham-Reichert}.
Its variational derivative is a constant multiple of the obstruction field $\mathcal{H}$,
which is an obstruction to the existence of $Y$ smooth up to the boundary;
see \cref{thm:asymptotically-minimal-submanifold,thm:first-variation-of-Graham-Witten-invariant}.
They have also shown that
any minimal submanifold in an Einstein manifold
is a critical point of the Graham-Witten energy~\cite[Proposition 4.5]{Graham-Reichert}.

In this paper,
we will consider the second variation of the Graham-Witten energy
at minimal submanifolds,
or more generally minimal immersions in an Einstein manifold.
Let $(M, g)$ be an $n$-dimensional Einstein manifold with Einstein constant $2 \lambda (n - 1)$,
and $f$ be a minimal immersion from a closed manifold $\Sigma$ of dimension $d = 2 k$ to $M$.
Consider a smooth map $F \colon (- \delta, \delta) \times \Sigma \to M$
such that $F_{t} = F(t, \cdot)$ is an immersion for each $t \in (- \delta, \delta)$ and $F_{0} = f$.
Assume that the variation field $V = F_{*} (\partial_{t}) |_{t = 0}$ is normal.
Then the second variation of the Graham-Witten energy at $f$ is written
in terms of $k$, $\lambda$, and the Jacobi operator $J$.

\begin{theorem} \label{thm:second-variation-of-Graham-Witten-invariant}
	The Graham-Witten energy $\mathcal{E}_{t}$ of $F_{t}$ satisfies
	\begin{align}
		\left. \frac{d^{2} \mathcal{E}_{t}}{d t^{2}} \right|_{t = 0}
		= \frac{1}{2 k} \int_{\Sigma} g\left( \mathcal{J} V, V \right) d a_{\Sigma},
	\end{align}
	where $d a_{\Sigma}$ is the area density on $\Sigma$ with respect to $g$,
	and the operator $\mathcal{J}$ is given by
	\begin{equation}
		\mathcal{J}
		= \prod_{l = 1}^{k + 1} [J + 2 \lambda (k + l)(k - l + 1)].
	\end{equation}
\end{theorem}

Note that an analogous result for the total $Q$-curvature
has been obtained by Matsumoto~\cite[Theorem 0.2]{Matsumoto13},
and Guillarmou, Moroianu, and Schlenker~\cite[Corollary 7.7]{Guillarmou-Moroianu-Schlenker}.

As an application,
we will prove that
any even-dimensional totally geodesic sphere in the unit sphere
is a critical point of the Graham-Witten energy with non-negative second variation.

\begin{theorem} \label{thm:second-variation-at-the-sphere}
	The second variation of $\mathcal{E}$
	at a totally geodesic sphere $S^{2k}$ in the unit sphere $S^{2k + m}(1)$ is non-negative,
	and positive in directions transverse to the orbit of the conformal group of $S^{2k + m}(1)$.
\end{theorem}

This result is a higher dimensional and codimensional generalization
of \cite[Proposition 1.2]{Graham-Reichert}.
We will also study some variational properties of the Graham-Witten energy
at Clifford hypersurfaces in the unit sphere
(\cref{prop:second-variation-at-Clifford-hypersurface}).

This paper is organized as follows.
In \cref{section:preliminaries},
we introduce some basic notions related to immersions to Riemannian manifolds.
\cref{section:Graham-Witten-invariant} presents
a definition of the Graham-Witten energy and some known results.
In \cref{section:boundary-value-problem-for-the-Jacobi-operator},
we discuss a boundary value problem of the Jacobi operator,
which plays an important role for the proof of \cref{thm:second-variation-of-Graham-Witten-invariant}.
\cref{section:proof-of-theorem} is devoted to
the proof of \cref{thm:second-variation-of-Graham-Witten-invariant}.
In \cref{section:applications},
we study the second variation of the Graham-Witten energy at totally geodesic spheres and Clifford hypersurfaces in the unit sphere
as applications of \cref{thm:second-variation-of-Graham-Witten-invariant}.

\section*{Acknowledgements}
The author is grateful to his supervisor Kengo Hirachi
for various helpful suggestions.
This work was supported by JSPS Research Fellowship for Young Scientists,
JSPS KAKENHI Grant Number JP16J04653,
and the Program for Leading Graduate Schools, MEXT, Japan.

\section{Preliminaries} \label{section:preliminaries}


Let $(M, g)$ be a Riemannian manifold.
Denote by $\nabla$ the Levi-Civita connection with respect to $g$,
and by $R$ the curvature of $g$;
our convention of $R$ is
$R(X, Y) = \nabla_{X} \nabla_{Y} - \nabla_{Y} \nabla_{X} - \nabla_{[X, Y]}$.
Let $\Sigma$ be a manifold of dimension $d \geq 1$
and $f \colon \Sigma \to M$ be an immersion.
The pullback $f^{*} T M$
has the orthogonal decomposition $f^{*} T M = T \Sigma \oplus N \Sigma$;
the bundle $N \Sigma$ is called the \emph{normal bundle} of $\Sigma$.
For a section $X$ of $f^{*} T M$,
we denote by $X^{\top}$ (resp.\ $X^{\bot}$)
the component of $X$ in $T \Sigma$ (resp.\ $N \Sigma$).
The \emph{second fundamental form} $\RN{2}$ is the vector bundle morphism
$S^{2} T \Sigma \to N \Sigma$ given by $\RN{2}(X, Y) = ((f^{*} \nabla)_{X} Y)^{\bot}$,
and the \emph{mean curvature} $H$ is defined by $H = d^{-1} \tr_{g} \RN{2} \in \Gamma(N \Sigma)$.
An immersion $f \colon \Sigma \to M$ is said to be \emph{totally geodesic} (resp.\ \emph{minimal})
 if $\RN{2} \equiv 0$ (resp.\ $H \equiv 0$).
The \emph{shape operator} $S \colon T \Sigma \otimes N \Sigma \to T \Sigma$
is defined by $S(X, \nu) = - ((f^{*} \nabla)_{X} \nu)^{\top}$.
The second fundamental form and the shape operator satisfies the following relation:
\begin{equation}
	g(\RN{2}(X, Y), \nu)
	= g(Y, S(X, \nu)).
\end{equation}

Next,
consider a one-parameter family of immersions.
Let $F \colon (- \delta, \delta) \times \Sigma \to M$ be
a smooth map such that $F_{t} = F(t, \cdot)$ is an immersion
for each $t \in (- \delta, \delta)$
and $F_{0} = f$.
We define $V$ to be $F_{*} (\partial_{t}) |_{t = 0}$,
called the \emph{variation field}.
The manifold $(- \delta, \delta) \times \Sigma$
has the canonical foliation $\mathcal{F} = (\{ t \} \times \Sigma)_{t \in (- \delta, \delta)}$,
and write $T \mathcal{F}$ for the vector bundle on $(- \delta, \delta) \times \Sigma$
whose fiber is the tangent space of each leaf.
The pullback $F^{*} T M$ splits into the orthogonal direct sum
$F^{*} T M = T \mathcal{F} \oplus N \mathcal{F}$;
note that the restriction of $N \mathcal{F}$ to a leaf $\{ t \} \times \Sigma$
coincides with the normal bundle with respect to $F_{t}$.
The pullback $F^{*} \nabla$ induces a connection on $N \mathcal{F}$,
denoted by $D$.
The mean curvature $H_{t}$ with respect to $F_{t}$
depends smoothly on $t$,
and defines a smooth section of $N \mathcal{F}$.
Then we have
\begin{equation}
	D_{\partial_{t}} H_{t} |_{t = 0}
	= - d^{-1} (\Delta^{\bot} - \mathcal{R} - \mathcal{B}) V^{\bot} + D_{V^{\top}} H,
\end{equation}
where $\Delta^{\bot}$ is the normal bundle Laplacian
(our convention of the Laplacian is chosen so that $\Delta^{\bot}$ is non-negative),
and $\mathcal{R}$ and $\mathcal{B}$ are the sections of $\End(N \Sigma)$
given by
\begin{align}
	\mathcal{R} \nu
	&= \tr_{g} [T \Sigma \otimes T \Sigma \to N \Sigma ; (X, Y) \mapsto (R(\nu, X) Y)^{\bot}], \\
	\mathcal{B} \nu
	&= \tr_{g} [T \Sigma \otimes T \Sigma \to N \Sigma ; (X, Y) \mapsto \RN{2}(S(X, \nu), Y)];
\end{align}
see~\cite[Section 2]{Weiner78} for a proof.
Set $J = \Delta^{\bot} - \mathcal{R} - \mathcal{B}$,
called the \emph{Jacobi operator} or the \emph{stability operator}.

\section{Graham-Witten energy} \label{section:Graham-Witten-invariant}

Let $(M, C)$ be a conformal manifold of dimension $n \geq 3$,
and $g$ be a representative of $C$.
Since we consider only a formal theory near the boundary,
we take $\overline{M^{+}} = \clop{0}{\epsilon_{0}}_{r} \times M$
as $\overline{X}$ in the introduction.
There exists a metric $g^{+}$ on $M^{+} = (0, \epsilon_{0}) \times M$
of the form
\begin{equation}
	g^{+} = \frac{d r^{2} + g_{r}}{r^{2}}
\end{equation}
satisfying
\begin{equation}
	\Ric(g^{+}) + n g^{+} = O(r^{n - 2}),
\end{equation}
where $g_{r}$ is a one-parameter family of metrics on $M$ with $g_{0} = g$.
These conditions determine the coefficients of the Taylor expansion $g_{r}$ in $r$ modulo $O(r^{n})$.
We call such a metric $g^{+}$
an \emph{asymptotically hyperbolic Einstein metric}.

Consider an immersion $f$ from a $2 k$-dimensional closed manifold $\Sigma$ to $M$.
Then there exists an ``asymptotically minimal'' immersion
$f^{+} \colon \Sigma^{+} = (0, \epsilon_{0}) \times \Sigma \to M^{+}$
whose boundary value coincides with $f$.

\begin{theorem}[{\cite[Theorem 3.1]{Graham-Reichert}}] \label{thm:asymptotically-minimal-submanifold}
	Let $\Exp$ be the exponential map on $M$ with respect to $g$.
	There exists a one-parameter family $u^{r} \in \Gamma(N \Sigma)$ of the form
	\begin{equation}
		u^{r}
		= \sum_{l = 1}^{k} u^{(l)} r^{2 l} + \mathcal{H} r^{2 k + 2} \log r
	\end{equation}
	such that
	\begin{equation}
		f^{+} \colon \Sigma^{+} \to M^{+}
		; (r, p) \mapsto (r, \Exp u^{r}(p))
	\end{equation}
	is an immersion,
	and the mean curvature $H^{+}$ of $f^{+}$
	satisfies $H^{+} = O(r^{2 k + 3} \log r)$.
	The coefficients $u^{(l)}$ and $\mathcal{H}$ are independent
	of the higher order terms in $g_{r}$.
\end{theorem}

Consider the area density $d a_{\Sigma^{+}}$ on $\Sigma^{+}$
with respect to the metric induced by $g^{+}$.
Since the area of $\Sigma^{+}$ must be infinite,
we apply a renormalization procedure.
Take $0 < \epsilon < \epsilon_{0}$,
and consider the area $\Area((\epsilon, \epsilon_{0}) \times \Sigma)$
of $(\epsilon, \epsilon_{0}) \times \Sigma$,
which is finite.
As $\epsilon \to + 0$,
it has the following asymptotic expansion:
\begin{equation}
	\Area((\epsilon, \epsilon_{0}) \times \Sigma)
	= \sum_{l = 0}^{k - 1} b_{l} \epsilon^{- 2 k + 2 l} + K \log (1 / \epsilon) + O(1).
\end{equation}
Moreover,
the coefficient $K$ of $\log (1 / \epsilon)$ is independent of the choice of a representative $g$ of $C$,
and gives a conformal invariant of the immersion $f$.
The \emph{Graham-Witten energy} $\mathcal{E}$ of $f$ is defined as
\begin{equation}
	K
	= 2 c_{k} \mathcal{E}, \qquad
	c_{k} = (- 1)^{k} [2^{2 k} (k - 1)! k!]^{- 1}.
\end{equation}
The first variation of the Graham-Witten energy $\mathcal{E}$
has been computed by Graham and Reichert.

\begin{theorem}[{\cite[Theorem 4.3]{Graham-Reichert}}] \label{thm:first-variation-of-Graham-Witten-invariant}
	Let $F \colon (- \delta, \delta) \times \Sigma \to M$
	be a one-parameter family of immersions such that $F_{0} = f$,
	and $V$ be its variation field.
	Then the Graham-Witten energy $\mathcal{E}_{t}$ of $F_{t}$ satisfies
	\begin{equation}
		\left. \frac{d \mathcal{E}_{t}}{d t} \right|_{t = 0}
		= - \frac{1}{4 k c_{k + 1}} \int_{\Sigma} g(\mathcal{H}, V) d a_{\Sigma},
	\end{equation}
	where $d a_{\Sigma}$ is the area density on $\Sigma$ with respect to $g$.
	In particular if $\mathcal{H} \equiv 0$,
	the first variation vanishes.
\end{theorem}

From now on,
we always assume that $C$ contains an Einstein metric $g$ with Einstein constant $2 \lambda (n - 1)$.
Then,
for sufficiently small $\epsilon_{0} > 0$,
the metric
\begin{equation}
	g^{+}
	= \frac{d r^{2} + \left( 1 - \lambda r^{2} / 2 \right)^{2} g}{r^{2}}
\end{equation}
on $M^{+}$ satisfies the Einstein equation $\Ric(g^{+}) = - n g^{+}$.
In what follows,
we use this $g^{+}$ as an asymptotically hyperbolic Einstein metric on $M^{+}$.

We also assume that
$f$ is minimal with respect to $g$.
Then the immersion $f^{+} \colon \Sigma^{+} \to M^{+}$
defined by $f^{+}(r, p) = (r, f(p))$
is also a minimal immersion.
This implies that
we can choose $u^{r}$ in \cref{thm:asymptotically-minimal-submanifold}
as $u^{r} \equiv 0$.
Hence the area of $(\epsilon, \epsilon_{0}) \times \Sigma$ can be computed as follows:
\begin{align}
	\Area((\epsilon, \epsilon_{0}) \times \Sigma)
	&= \int_{\epsilon}^{\epsilon_{0}} \int_{\Sigma} r^{- 2 k -1} (1 - \lambda r^{2} / 2)^{2k} dr d a_{\Sigma} \\
	&= \sum_{l = 0}^{k -1} \frac{1}{2(k - l)} \binom{2 k}{l} \left(- \frac{\lambda}{2} \right)^{l} \Area(\Sigma) \epsilon^{- 2 (k - l)} \\
	&\quad + \binom{2 k}{k} \left(- \frac{\lambda}{2} \right)^{k} \Area(\Sigma) \log(1 / \epsilon) + O(1),
\end{align}
where $\Area(\Sigma)$ is the area of $\Sigma$ with respect to $g$.
This gives an explicit formula of $\mathcal{E}$,
which may be used implicitly in \cite{Graham-Reichert}:

\begin{proposition}
	The Graham-Witten energy $\mathcal{E}$ of $f$ is given by
	\begin{equation}
		\mathcal{E}
		= 2^{k} (2 k - 1)! \lambda^{k} \Area(\Sigma).
	\end{equation}
\end{proposition}

Moreover,
the obstruction field $\mathcal{H}$ vanishes identically,
and consequently the first variation of $\mathcal{E}$ at $f$ vanishes~\cite[Proposition 4.5]{Graham-Reichert}.

\section{Boundary value problem for the Jacobi operator} \label{section:boundary-value-problem-for-the-Jacobi-operator}

In this section,
we consider a boundary value problem for the Jacobi operator on $\Sigma^{+}$,
which may be of independent interest
as well as will play a crucial role in the proof of \cref{thm:second-variation-of-Graham-Witten-invariant} in the next section.
A normal vector field on $\overline{\Sigma^{+}} = \clop{0}{\epsilon_{0}} \times \Sigma$
with respect to $\overline{g} = d r^{2} + (1 - \lambda r^{2} / 2)^{2} g$
can be identified with a family of normal vector fields on $\Sigma$
depending smoothly on $r \in \clop{0}{\epsilon_{0}}$.
Under this identification,
we prove

\begin{theorem} \label{thm:scattering-for-Jacobi-operator}
	For $\nu \in \Gamma(N \Sigma)$,
	there exist $\nu^{(1)} \in \Gamma(N \overline{\Sigma^{+}})$ with $\nu^{(1)} |_{r = 0} = \nu$
	and $\nu^{(2)} \in \Gamma(N \overline{\Sigma^{+}})$ such that
	\begin{equation}
		\nu^{+} = \nu^{(1)} + \nu^{(2)} r^{2 k + 2} \log r \in \Gamma(N \Sigma^{+})
	\end{equation}
	satisfies the equation
	\begin{equation}
		J^{+} \nu^{+} = O(r^{2 k + 3} \log r),
	\end{equation}
	where $J^{+}$ is the Jacobi operator for $f^{+}$.
	Moreover,
	$\nu^{(1)}$ modulo $O(r^{2 k + 2})$ and $\nu^{(2)} |_{r = 0}$ depend only on $\nu$,
	and
	\begin{equation}
		\nu^{(2)} |_{r = 0}
		= - 2 c_{k + 1} \mathcal{J} \nu,
	\end{equation}
	where the operator $\mathcal{J}$ is as in \cref{thm:second-variation-of-Graham-Witten-invariant}.
\end{theorem}

Our proof is similar in spirit to the proof given by Fefferman and Graham~\cite[Proposition 7.9]{Fefferman-Graham12}
of an analogous formula for GJMS operators on Einstein manifolds.
We first show the existence and uniqueness of $\nu^{(1)}$ and $\nu^{(2)}$.
Let $h(r)$ be a smooth function in $r$,
and $\nu$ be a section of $N \Sigma$.
Then a straightforward calculation shows that
\begin{align}
	&(1 - \lambda r^{2} / 2)^{2} J^{+} (h \cdot \nu) \\
	&= r^{2} h \cdot J \nu - (1 - \lambda r^{2} / 2)^{2} (r \partial_{r})^{2} h \cdot \nu
		+ 2(k + 1)(1 - \lambda^{2} r^{4} / 4) (r \partial_{r}) h \cdot \nu.
\end{align}
Thus we can determine $\nu^{(1)}$ and $\nu^{(2)}$ inductively
as in the proof of \cite[Proposition 4.2]{Graham-Zworski03}.
Moreover,
from the construction,
we obtain a polynomial $p_{k} \in \mathbb{R}[x, \lambda]$
such that its coefficients depend only on $k$,
$p_{k}(x, \lambda) = \lambda^{k + 1} p_{k}(x / \lambda, 1)$ for $\lambda \neq 0$,
and $\nu^{(2)} |_{r = 0} = p_{k}(J, \lambda) \nu$.

Therefore,
it suffices to obtain an explicit formula of $p_{k}$.
The homogeneity of $p_{k}$ implies that
$p_{k}$ is determined by $p_{k}(x, 1 / 2)$.
Take as $M$ the unit sphere $S^{2 k + 1}(1)$ in $\mathbb{R}^{2 k + 2}$,
and as $\Sigma$ the totally geodesic sphere in $S^{2 k + 1}(1)$ defined by $\{ x^{2 k + 2} = 0 \}$.
Set $\nu_{1} = \partial / \partial x^{2 k + 2}$,
which is a parallel normal vector field on $\Sigma$ with unit length.
Then $\nu_{1}^{+} =  (1 - \lambda r^{2} / 2)^{-1} r \nu_{1}$
gives a parallel normal vector field on $\Sigma^{+}$
with unit length.
For a section $\nu^{+} = \phi^{+} \nu_{1}^{+}$,
\begin{equation}
	J^{+} \nu^{+}
	= [(\Delta^{+} + 2 k + 1) \phi^{+}] \nu_{1}^{+},
\end{equation}
where $\Delta^{+}$ is the scalar Laplacian with respect to $g^{+} |_{\Sigma^{+}}$,
which is also an asymptotically hyperbolic Einstein metric on $\Sigma^{+}$.
By \cite[Propositions 4.2 and 4.3]{Graham-Zworski03},
for $\phi \in C^{\infty}(\Sigma)$,
there exists $\phi^{(1)}, \phi^{(2)} \in C^{\infty}(\overline{\Sigma^{+}})$ such that
\begin{equation}
	(\Delta^{+} + 2 k + 1) (\phi^{(1)} r^{-1} + \phi^{(2)} r^{2 k + 1} \log r) = O(r^{2 k + 2} \log r),
	\quad \phi^{(1)} |_{r = 0} = \phi.
\end{equation}
Moreover,
$\phi^{(2)} |_{r = 0} = - 2 c_{k + 1} P_{k + 1} \phi$,
where
\begin{equation}
	P_{k + 1}
	= \prod_{l = 1}^{k + 1} [\Delta + (k + l - 1)(k - l)]
\end{equation}
and $\Delta$ is the scalar Laplacian on $S^{2 k}(1)$;
see \cite[Theorem 1.2]{Gover06} and \cite[Proposition 7.9]{Fefferman-Graham12}.
Hence if we set
\begin{equation}
	\nu^{(1)}
	= (1 - \lambda r^{2} / 2)^{-1} \phi^{(1)} \nu_{1}, \quad
	\nu^{(2)}
	= (1 - \lambda r^{2} / 2)^{-1} \phi^{(2)} \nu_{1},
\end{equation}
$\nu^{+} = \nu^{(1)} + \nu^{(2)} r^{2 k + 2} \log r$ satisfies
$J^{+} \nu^{+} = O(r^{2 k + 3} \log r)$ and $\nu^{(1)} |_{r = 0} = \phi \nu_{1}$.
Moreover,
$\nu^{(2)} |_{r = 0}$ corresponding to $\nu = \phi \nu_{1}$
is given by
\begin{equation}
	\nu^{(2)} |_{r = 0}
	= - 2 c_{k + 1} \left\{ \prod_{l = 1}^{k + 1} [\Delta + (k + l - 1)(k - l)] \phi \right\} \nu_{1}.
\end{equation}
From $J(\phi \nu_{1}) = [(\Delta - 2 k) \phi] \nu_{1}$,
we obtain
\begin{equation}
	\nu^{(2)} |_{r = 0}
	= - 2 c_{k + 1} \left\{ \prod_{l = 1}^{k + 1} [J + (k + l)(k - l + 1)] \right\} (\phi \cdot \nu_{1}).
\end{equation}
Therefore,
we have
\begin{equation}
	p_{k}(x, \lambda)
	= - 2 c_{k + 1} \prod_{l = 1}^{k + 1} [x + 2 \lambda (k + l)(k - l + 1)].
\end{equation}
This completes the proof of \cref{thm:scattering-for-Jacobi-operator}.

\section{Proof of \cref{thm:second-variation-of-Graham-Witten-invariant}} \label{section:proof-of-theorem}

In this section,
we give a proof of \cref{thm:second-variation-of-Graham-Witten-invariant}.
Consider a one-parameter family $F \colon (- \delta, \delta) \times \Sigma \to M$ of immersions
such that $F_{0} = f$ and the variation field $V = F_{*} (\partial_{t}) |_{t = 0}$ is normal.
According to \cref{thm:first-variation-of-Graham-Witten-invariant},
the second variation is written as
\begin{equation}
	\left. \frac{d^{2} \mathcal{E}_{t}}{d t^{2}} \right|_{t = 0}
	= - \frac{1}{4 k c_{k + 1}} \int_{\Sigma} g(D_{\partial_{t}} \mathcal{H}_{t} |_{t = 0}, V) d a_{\Sigma}.
\end{equation}
Hence
it is sufficient to compute $D_{\partial_{t}} \mathcal{H}_{t} |_{t = 0}$.

Let
\begin{equation}
	U^{r}_{t}
	= \sum_{l = 1}^{k} U^{(l)}_{t} r^{2 l} + \mathcal{H}_{t} r^{2 k + 2} \log r
\end{equation}
be a normal vector field as in \cref{thm:asymptotically-minimal-submanifold}
with respect to $F_{t}$ such that $U^{r}_{0} \equiv 0$.
Define a map $F^{+} \colon (- \delta, \delta) \times \Sigma^{+} \to M^{+}$ by
\begin{equation}
	F^{+}(t, r, p)
	= (r, \Exp U^{r}_{t} (p)),
\end{equation}
which gives a smooth family of immersions from $\Sigma^{+}$ to $M^{+}$.
Then the variation field $V^{+} = (F^{+})_{*}(\partial_{t}) |_{t = 0}$ of $F^{+}$ is given by
\begin{equation}
	V^{+}
	= V + \sum_{l = 1}^{k} D_{\partial_{t}} U^{(l)}_{t} |_{t = 0} r^{2 l}
	+ D_{\partial_{t}} \mathcal{H}_{t} |_{t = 0} r^{2 k + 2} \log r,
\end{equation}
which is normal.
Hence
the mean curvature $H^{+}_{t}$ of $F^{+}_{t}$ satisfies
\begin{equation}
	(D^{+})_{\partial_{t}} H^{+}_{t} |_{t = 0}
	= - (2 k + 1)^{-1} J^{+} V^{+}.
\end{equation}
From $H^{+}_{t} = O(r^{2 k + 3} \log r)$,
it follows that $J^{+} V^{+} = O(r^{2 k + 3} \log r)$.
Therefore,
\cref{thm:scattering-for-Jacobi-operator} implies
$D_{\partial_{t}} \mathcal{H}_{t} |_{t = 0} = - 2 c_{k + 1} \mathcal{J} V$,
which proves \cref{thm:second-variation-of-Graham-Witten-invariant}.

\section{Applications} \label{section:applications}

In this section,
we study the second variation of the Graham-Witten energy
at totally geodesic spheres and Clifford hypersurfaces in the unit sphere.

We first recall some properties of the unit sphere $S^{n} (1)$.
Let $\mathfrak{conf}(S^{n}(1))$ be the space of conformal vector fields on $S^{n}(1)$.
This space is isomorphic to
\begin{equation}
	\mathfrak{so}(n + 1, 1)
	=
	\Set{
		\begin{pmatrix}
			0 & ^{t} b \\
			b & A
		\end{pmatrix}
		\in \mathfrak{gl}(n + 2, \mathbb{R})
		|
		A \in \mathfrak{so}(n + 1),
		b \in \mathbb{R}^{n + 1}
	}
\end{equation}
by the following map:
\begin{align}
	\mathfrak{so}(n + 1, 1) &\to \mathfrak{conf}(S^{n}(1)) \\
	\begin{pmatrix}
		0 & ^{t} b \\
		b & A
	\end{pmatrix}
	&\mapsto
	[x \mapsto A x - \Braket{b, x} x + b \in (\mathbb{R} x)^{\bot} \cong T_{x} S^{n}(1)],
\end{align}
where $\Braket{\cdot, \cdot}$ is the standard inner product on $\mathbb{R}^{n + 1}$.
In particular,
the dimension of $\mathfrak{conf}(S^{n}(1))$ is $(n + 1)(n + 2) / 2$.
The space $\mathfrak{so}(n + 1, 1)$ has the direct sum decomposition
$\mathfrak{so}(n + 1, 1) \cong \mathfrak{so}(n + 1) \oplus \mathbb{R}^{n + 1}$;
the image of $\mathfrak{so}(n + 1)$ by the above map
is the space $\mathfrak{isom}(S^{n} (1))$ of Killing vector fields on $S^{n}(1)$,
and that of $\mathbb{R}^{n + 1}$ coincides with the space $\mathfrak{const}(S^{n}(1))$
of tangential projections onto $S^{n} (1)$ of constant vector fields on $\mathbb{R}^{n + 1}$.

Let $\Sigma \subset S^{n} (1)$ be a minimal submanifold of dimension $d$,
and consider the image $\mathscr{C}$ of $\mathfrak{conf}(S^{n} (1))$
by the canonical restriction-projection map
\begin{equation}
	\Gamma(T S^{n} (1)) \to \Gamma(T S^{n} (1) |_{\Sigma}) \to \Gamma(N \Sigma).
\end{equation}
The image of $\mathfrak{isom}(S^{n} (1))$ is contained in $\ker J$,
while that of $\mathfrak{const}(S^{n}(1))$ is annihilated by $J + d$;
see \cite[Section 5.1]{Simons68}.
In particular,
we have $\mathscr{C} \subset \ker J (J + d)$.

Now we consider a totally geodesic sphere $\Sigma = S^{2 k} \subset M = S^{2 k + m}(1)$.
From the conformal invariance of the Graham-Witten energy,
it follows that $\mathscr{C}$ is contained in the kernel of $\mathcal{J}$.
The kernel of the map $\mathfrak{so}(2 k + m + 1, 1) \cong \mathfrak{conf}(S^{2 k + m}(1)) \to \Gamma(N S^{2 k})$
is the space
\begin{equation}
	\Set{
		\begin{pmatrix}
			0 & ^{t} b_{1} & 0 \\
			b_{1} & A_{1} & 0 \\
			0 & 0 & A_{2}
		\end{pmatrix}
		|
		A_{1} \in \mathfrak{so}(2 k + 1),
		A_{2} \in \mathfrak{so}(m),
		b_{1} \in \mathbb{R}^{2 k + 1}
	},
\end{equation}
whose dimension is $(k + 1)(2 k + 1) + m(m - 1) / 2$.
Thus
we obtain $\dim \mathscr{C} = 2 (k + 1) m$.

\begin{proof}[Proof of \cref{thm:second-variation-at-the-sphere}]
	Set $\nu_{i} = \partial / \partial x^{2 k + i + 1} (i = 1, \dots , m)$,
	which give a parallel orthonormal frame of $N S^{2 k}$.
	Any smooth section $V$ of $N S^{2 k}$
	is of the form $\sum_{i = 1}^{m} \phi_{i} \nu_{i}$,
	$\phi_{i} \in C^{\infty}(S^{2 k})$.
	Since $J (\phi_{i} \nu_{i}) = [(\Delta - 2 k) \phi_{i}] \nu_{i}$,
	we obtain
	\begin{equation}
		\mathcal{J} V
		= \sum_{i = 1}^{m} \left\{\prod_{l = 1}^{k + 1} [\Delta + (k + l - 1)(k - l)] \phi_{i} \right\} \nu_{i}.
	\end{equation}
	The eigenvalues of $\Delta$ on $S^{2 k}(1)$
	are $j (2 k + j - 1)$,
	$j \in \mathbb{N}$,
	with multiplicity $\binom{2 k + j}{2 k} - \binom{2 k + j - 2}{2 k}$.
	This fact gives that
	$\mathcal{J}$ is non-negative
	and its kernel $\ker \mathcal{J}$ has dimension $2 (k + 1) m$.
	Since $\mathscr{C} \subset \ker \mathcal{J}$
	and $\dim \mathscr{C} = 2 (k + 1) m$,
	we have $\mathscr{C} = \ker \mathcal{J}$.
	This proves the positivity of the second variation in directions
	transverse to the orbit of the conformal group.
\end{proof}

We next turn to studying variational properties of the Graham-Witten energy at Clifford hypersurfaces.
Let $d_{1}$ and $d_{2}$ be positive integers with even $d_{1} + d_{2}$,
denoted by $2 k$,
and set $r_{i} = \sqrt{d_{i} / (2 k)}$.
Consider the sphere $S^{d_{i}} (r_{i})$ in $\mathbb{R}^{d_{i} + 1}$ with radius $r_{i}$.
The product $\Sigma = S^{d_{1}} (r_{1}) \times S^{d_{2}} (r_{2})$ in $\mathbb{R}^{2 k + 2}$
defines a minimal hypersurface in $M = S^{2 k + 1} (1)$,
called a \emph{Clifford hypersurface}.

As in the previous case,
consider the image $\mathscr{C}$ of $\mathfrak{conf}(S^{2 k + 1} (1))$ by the restriction-projection map.
The kernel of $\mathfrak{so}(2 k + 2, 1) \cong \mathfrak{conf}(S^{2 k + 1} (1)) \to \Gamma(N \Sigma)$
is the space
\begin{equation}
	\Set{
		\begin{pmatrix}
			0 &  0 & 0 \\
			0 & A_{1} & 0  \\
			0 & 0 & A_{2}  
		\end{pmatrix}
		|
		A_{i} \in \mathfrak{so}(d_{i} + 1)
	}.
\end{equation}
Hence the dimension of $\mathscr{C}$ is given by
\begin{align}
	\dim \mathscr{C}
	&= (2 k + 2) (2 k + 3) / 2
	- d_{1} (d_{1} + 1) / 2 - d_{2} (d_{2} + 1) / 2 \\
	&= (d_{1} + d_{2} + 2) + (d_{1} + 1)(d_{2} + 1).
\end{align}

To study the spectrum of the Jacobi operator on $\Sigma$,
we first construct a parallel frame of $N \Sigma$.
The Euler vector field $E_{i}$ on $\mathbb{R}^{d_{i} + 1}$
gives a parallel frame of the normal bundle $N S^{d_{i}}(r_{i})$ with length $r_{i}$.
Then the vector field $\mu = d_{2} E_{1} - d_{1} E_{2}$
defines a parallel frame of $N \Sigma$.
A calculation shows that $\mathcal{R} = \mathcal{B} = 2 k$,
and $J (\psi \mu) = [(\Delta - 4 k) \psi] \mu$,
where $\Delta$ is the scalar Laplacian on $\Sigma$.
Thus we have
\begin{equation}
	\mathcal{J}\left(\psi \mu \right)
	= \left\{\prod_{l = 1}^{k + 1}
	[\Delta - 4 k + (k + l) (k - l + 1)] \psi \right\} \mu,
\end{equation}
and consequently it is enough to study the spectrum of
\begin{equation}
	L
	= \prod_{l = 1}^{k + 1} [\Delta - 4 k + (k + l) (k - l + 1)].
\end{equation}

\begin{proposition} \label{prop:second-variation-at-Clifford-hypersurface}
	The second variation of $\mathcal{E}$ at $\Sigma$ is positive
	in all directions orthogonal to $\mu$ and to the orbit of the conformal group of $S^{2 k + 1} (1)$.
	The second variation in the direction $\mu$
	is positive if $k = 1$, zero if $k = 3$, and negative if $k = 2$ or $k \geq 4$.
\end{proposition}

Note that this result has been already obtained
by Weiner~\cite[Proposition 3.1]{Weiner78}, and Graham and Reichert~\cite[Proposition 6.5]{Graham-Reichert}
for the case of $d_{1} = d_{2} = 1$ and $d_{1} = d_{2} = 2$, respectively.

\begin{proof}[Proof of \cref{prop:second-variation-at-Clifford-hypersurface}]
	If $\psi$ is an eigenfunction of $\Delta$ with eigenvalue $\lambda$,
	then it is also an eigenfunction of $L$ with eigenvalue
	\begin{equation} \label{eq:eigenvalue-of-J}
		\prod_{l = 1}^{k + 1} [\lambda - 4 k + (k + l) (k - l + 1)].
	\end{equation}
	In particular if $\lambda > 4 k$,
	this value is positive.
	Thus,
	we consider only the eigenvalues of $\Delta$ smaller than or equal to $4 k$.
	
	The eigenvalues of the scalar Laplacian on $S^{d_{i}} (r_{i})$ are $r_{i}^{- 2} j (d_{i} + j - 1)$,
	$j \in \mathbb{N}$,
	with multiplicity $\binom{d_{i} + j}{d_{i}} - \binom{d_{i} + j - 2}{d_{i}}$.
	Hence
	the eigenvalues of $\Delta$ at most $4 k$
	are $0$, $2 k$, $4 k$
	with multiplicity $1$, $d_{1} + d_{2} + 2$, $(d_{1} + 1)(d_{2} + 1)$, respectively.
	In particular,
	the dimension of $\ker (\Delta - 2 k)(\Delta  - 4 k)$
	coincides with that of $\mathscr{C}' = \Set{ \psi \in C^{\infty}(\Sigma) \mid \psi \mu \in \mathscr{C} }$.
	On the other hand,
	as we already noted, 
	$\mathscr{C} \subset \ker J (J + 2 k)$,
	which is equivalent to $\mathscr{C}' \subset \ker (\Delta - 2 k) (\Delta - 4 k)$.
	Therefore,
	we have $\mathscr{C}' = \ker (\Delta - 2 k) (\Delta - 4 k) \subset \ker L$.
	For $\lambda = 0$,
	or constant $\psi$,
	the sign of \cref{eq:eigenvalue-of-J} depends on $k$;
	it is positive if $k = 1$, zero if $k = 3$, and negative if $k = 2$ or $k \geq 4$.
	This completes the proof.
\end{proof}

\end{document}